\documentclass[a4paper,12pt]{article}

\usepackage{color,makeidx,bm,latexsym,amscd,amsmath,amssymb,enumerate,stmaryrd,amsthm,float,graphicx,psfrag,geometry,comment}

\usepackage[all]{xy}
\usepackage{bm}

\usepackage{tikz,epic,eso-pic}

\tikzset{v/.style={
  circle, draw, inner sep=2pt, minimum size=6pt, fill=white}}

\theoremstyle{plain}
\newtheorem{theorem}{Theorem}[section]
\newtheorem{corollary}[theorem]{Corollary}

\newtheorem{conjecture}[theorem]{Conjecture}

\theoremstyle{definition}
\newtheorem{definition}[theorem]{Definition}
\newtheorem{remark}[theorem]{Remark}
\newtheorem{example}[theorem]{Example}
\newtheorem{proposition}[theorem]{Proposition}

\def\qed{\hfill $\Box$}

\newcommand{\Img}{\operatorname{Im}}
\newcommand{\Ker}{\operatorname{Ker}}
\newcommand{\ch}{\operatorname{ch}}
\newcommand{\bch}{\operatorname{bch}}
\newcommand{\uch}{\operatorname{uch}}

\newcommand{\Sep}{\operatorname{Sep}}
\newcommand{\Hom}{\operatorname{Hom}}
\newcommand{\Gal}{\operatorname{Gal}}

\newcommand{\bC}{\mathbb{C}}

\newcommand{\bP}{\mathbb{P}}

\newcommand{\bR}{\mathbb{R}}

\newcommand{\bZ}{\mathbb{Z}}

\newcommand{\cA}{\mathcal{A}}

\newcommand{\cD}{\mathcal{D}}

\newcommand{\cF}{\mathcal{F}}

\newcommand{\cI}{\mathcal{I}}

\newcommand{\cL}{\mathcal{L}}

\newcommand{\cS}{\mathcal{S}}

\def\qed{\hfill $\Box$}

\title{$\mathbb{Z}$-local system cohomology of hyperplane arrangements and a Cohen-Dimca-Orlik type theorem}
\author{Sakumi Sugawara}
\date{\today}

\begin{document}
\maketitle

\begin{abstract}
Local system cohomology groups of the complements of hyperplane arrangements have played an important role in the theory of hypergeometric integrals, the topology of Milnor fibers and covering spaces. One of the important theorems is the vanishing theorem for generic $\bC$-local systems which goes back to Aomoto's work. Later, Cohen, Dimca, and Orlik proved a stronger version of the vanishing theorem. In this paper, we prove a Cohen-Dimca-Orlik type theorem for $\bZ$-local systems.
\end{abstract}

\section{Introduction}

A finite set $\cA =\{H_1 , \cdots , H_n\}$ of hyperplanes in a vector space $\bC^{\ell}$ is called a hyperplane arrangement. A hyperplane arrangement defines a poset of nonempty intersections of hyperplanes, called the intersection poset.
Hyperplane arrangements are studied from several viewpoints, for example, topology, combinatorics, and algebraic geometry (see \cite{orl-ter-arr} for details).
From topological aspects, one of the important questions is whether topology of the complement $M(\cA) = \bC^{\ell} \setminus \bigcup_{i=1}^{n} H_i$ is determined by combinatorial structure, that is, the information of intersection poset.
For example, Betti numbers, the cohomology ring, and local system cohomology for generic local systems are combinatorially determined. 
Recently, Yoshinaga proved that mod $2$ Betti numbers of the double covering are combinatorially determined (\cite{yos-dou}).
Whether characteristic varieties and Betti numbers of Milnor fibers are combinatorially determined are still open problems. 
See \cite{pap-suc} for a recent study on these problems.

Local system cohomology groups of the complements are well studied in the theory of hypergeometric integrals \cite{orl-ter-int}, and topology of Milnor fibers and covering spaces \cite{coh-suc-mil}.
There are some kinds of vanishing theorems (\cite{aom,koh-hom,lib-eig,cdo}), which say that under some generic conditions, the local system cohomology group concentrates in the top degree and vanishes in other degrees.
Cohen, Dimca, and Orlik proved the following theorem for rank one $\bC$-local systems.

\begin{theorem} \label{thm:cdo} {\rm (\cite{cdo})}
Suppose that a rank one $\bC$-local system $\cL$ satisfies CDO-condition. Then,
\begin{eqnarray*}
H^{k} (M(\cA), \cL) \cong
\left \{
\begin{array} {ll}
\bC^{|\chi (M(\cA))|} & (k= \ell ), \\
0 & (k \neq \ell),
\end{array}
\right.
\end{eqnarray*}
where $\chi (M(\cA))$ is the Euler characteristic of $M(\cA)$.
\end{theorem}

The definition of the CDO-condition is given in Definition \ref{def:cdo}. Since the Euler characteristic $\chi (M(\cA))$ is determined by the intersection poset, we have the following.

\begin{corollary}
If a rank one $\bC$-local system $\cL$ satisfies the CDO-condition, then the local system cohomology $H^{*}( M(\cA), \cL)$ is combinatorially determined. 
\end{corollary}

In this paper, we focus on rank one $\bZ$-local systems, that is, the local system whose coefficient modules are locally isomorphic to $\bZ$. A rank one $\bZ$-local system is determined by the characteristic homomorphism $\rho: H_1 (M(\cA),\bZ) \rightarrow \bZ^{\times} = \{\pm1\}$. 
We prove a Cohen-Dimca-Orlik type theorem for $\bZ$-local system cohomology for complexified real arrangements. 
We calculate $\bZ$-local system cohomology.
The main result of this paper is the following.

\begin{theorem} 
{\rm (Theorem 3.1)}
Let $\cA = \{H_1 , \cdots , H_n\}$ be a real hyperplane arrangement in $\bR^{\ell}$ and denote by $M(\cA) $ the complexified complement. Suppose that the rank one $\bZ$-local system $\cL$ satisfies CDO-condition. Then,
\begin{eqnarray*}
H^{k} (M(\cA), \cL) \cong
\left \{
\begin{array} {ll}
\bZ^{\beta_{\ell}} \oplus \bZ_2^{\beta_{\ell -1}} & (k= \ell ), \\
\bZ_2^{\beta_{k-1}} & (1 \leq k \leq \ell-1), \\
0 & (otherwise),
\end{array}
\right.
\end{eqnarray*}
where $\beta_{k} = |\sum_{i=0}^{k} (-1)^i \cdot b_i (M(\cA))|$ and $\bZ_2 = \bZ / 2\bZ$.
\end{theorem}

A rank one $\bZ$-local system $\cL$ defines a rank one $\bC$-local system $\cL \otimes_{\bZ} \bC$. Therefore, by tensoring with $\bC$, we obtain the following:
\begin{eqnarray*}
H^{k} (M(\cA), \cL \otimes_{\bZ} \bC) \cong
\left \{
\begin{array} {ll}
\bC^{\beta_{\ell}} & (k= \ell ), \\
0 & (k \neq \ell).
\end{array}
\right.
\end{eqnarray*}
Note that $\beta_{\ell} = |\sum_{i=0}^{\ell} (-1)^i \cdot b_i (M(\cA))| = |\chi (M(\cA))|$. This agrees with Theorem \ref{thm:cdo}. 

The proof is based on the minimality of $M(\cA)$ and the description of the twisted cochain complex by using the real figure of chambers.
These methods were developed in \cite{yos-lef,yos-ch,yos-loc,yos-mil,bai-yos} and applied to calculate monodromy eigenspaces of Milnor fibers and to prove vanishing results for the Aomoto complex.
In this paper, we calculate $\bZ$-local system cohomology by using these methods.

Since Betti numbers are combinatorially determined, the main result implies that $H^{*} (M(\cA), \cL)$ is combinatorially determined for complexified real arrangements. 
This leads us the following conjecture.

\begin{conjecture} \label{conjecture}
Let $\cA$ be a complex hyperplane arrangement in $\bC^{\ell}$ (not necessarily a complexified real arrangement).
If a rank one $\bZ$-local system $\cL$ on $M(\cA)$ satisfies the CDO-condition, then
\begin{eqnarray*}
H^{k} (M(\cA), \cL) \cong
\left \{
\begin{array} {ll}
\bZ^{\beta_{\ell}} \oplus \bZ_2^{\beta_{\ell -1}} & (k= \ell ), \\
\bZ_2^{\beta_{k-1}} & (1 \leq k \leq \ell-1), \\
0 & (otherwise),
\end{array}
\right.
\end{eqnarray*}
where $\beta_{k} = |\sum_{i=0}^{k} (-1)^i \cdot b_i (M(\cA))|$. Thus, the local system cohomology $H^{*} (M(\cA),\cL)$ is combinatorially determined.
\end{conjecture}

\begin{remark}
Recently, Liu, Maxim, and Wang wrote the preprint claiming that Conjecture \ref{conjecture} holds (\cite{lmw}). They use the theory of perverse sheaves to calculate the cohomology.
\end{remark}

This paper is organized as follows. In Section $2$, we review basic notations of real arrangements and the methods to calculate local system cohomology in terms of chamber basis. 
In Section $3$, we give the proof of the main theorem. 
In Section $4$, we give some examples of $\bZ$-local system cohomology, including the icosidodecahedral arrangement. 










\medskip

\noindent
{\bf Acknowledgement.}
The author would like to thank Masahiko Yoshinaga for helpful discussions and comments on the paper. This work was supported by JSPS KAKENHI Grant Number 22J20470.

\section{Local system cohomology on $M(\cA)$}



\subsection{Local system on arrangements}
Let $\cA = \{H_1, \cdots , H_n\}$ be a hyperplane arrangement in $\bR^{\ell}$. We denote by $M(\cA) = \bC^{\ell} \setminus \bigcup_{i=1}^{n} (H_i \otimes \bC)$ be the complexified complement. 
We identify the affine space $\bR^{\ell}$ with $\bR \bP^{\ell} \setminus \overline{H}_{\infty}$, and let $\overline{\cA} =\{\overline{H}_1 , \cdots , \overline{H}_n , \overline{H}_{\infty} \}$ be the projective closure of $\cA$, where $\overline{H}_{\infty}$ denotes the hyperplane at infinity. 
The poset consisting of nonempty intersections of hyperplanes in $\cA$ is denoted by $L(\cA)$. 
For each $0 \leq k \leq \ell$, denote by $L_{k} (\cA)$ the set of $k$-dimensional intersection of $\cA$.
For example, $L_{\ell} (\cA)= \{ \bR^{\ell} \}$ and $L_{\ell - 1} (\cA) = \cA$.
The intersection $X \in L(\cA)$ is called a {\it dense edge} if the localization $\cA_{X} = \{H \supset X \mid  H \in \cA \}$ is indecomposable, that is $\cA_{X}$ cannot be expressed as the direct sum of nonempty two subarrangements.
The intersection poset $L(\overline{\cA})$ and dense edges in the projective closure $\overline{\cA}$ are defined similarly.
The set of dense edges contained in $\overline{H}_{\infty}$ is denoted by $D_{\infty} (\cA)$, that is, 
\[
D_{\infty} (\cA) = \{\mbox{$X$ : dense edge} \mid X \subset \overline{H}_{\infty} \}.
\]

Let $R = \bC$ or $\bZ$ be the coefficient ring. A rank one $R$-coefficient local system $\cL$ on $M(\cA)$ is determined by the characteristic homomorphism $\rho : H_1 (M(\cA) , \bZ) \rightarrow R^{\times }$.
This is equivalent to determining a $n$-tuple $(q_1 , \cdots , q_n) \in (R^{\times})^{n} $ where $q_i = \rho (\gamma_i)$ is the local monodromy of the meridian $\gamma_i$ around the hyperplane $H_i$. 
Denote $q_{X} = \prod_{H_i \supset X} q_i$ for each intersection $X \in L(\cA)$, and $q_{\infty} = \prod_{i=1}^{n} q_i ^{-1} $ the monodromy around the hyperplane at infinity.


Now, let us recall the CDO-condition for a local system $\cL$.

\begin{definition} \label{def:cdo}
 Let $\cL$ be a rank one $R$-local system on $M(\cA)$. We call that $\cL$ {\it satisfies the CDO-condition} if $q_X \neq 1$ for each $X \in D_{\infty} (\cA)$.
\end{definition}



\subsection{Chambers on real arrangements}
Let $\cA = \{H_1, \cdots , H_n\}$ be a real arrangement in $\bR^{\ell}$. 
A connected component of the complement $\bR^{\ell} \setminus \bigcup_{i=1}^{n} H_i$ is called a chamber, and denote by $\ch(\cA)$ the set of all chambers.
Each chamber is either bounded or unbounded. 
Let $\bch(\cA)$ be the set of all bounded chambers and $\uch (\cA)$ be the set of unbounded chambers. 
It is clear that $\ch (\cA) = \bch (\cA) \sqcup \uch(\cA)$.
Remark that $C \in \uch (\cA)$ is equivalent to $\overline{C} \cap \overline{H}_{\infty} \neq \emptyset$, where $\overline{C}$ denotes the closure of $C$ in the projective space.

\begin{definition}
Let $C$, $C' \in \ch (\cA)$. The set $\Sep (C, C')$ of separating hyperplanes is defined as follows: 
\[
\Sep (C, C') = \{H \in \cA \mid \mbox{$H$ separates $C$ and $C'$} \}.
\]
\end{definition}

Let $\cF$ be a generic flag in $\bR^{\ell}$, that is, $\cF$ is a sequence of affine subspaces in $\bR^{\ell}$ 
\[
\cF : \emptyset = F^{-1} \subset F^{0} \subset \cdots \subset F^{\ell-1} \subset F^{\ell} = \bR^{\ell},
\]
where $F^{k}$ is a $k$-dimensional subspace such that $\dim (X \cap F^{k}) = k + \dim X - \ell$ for each $X \in L(\cA)$ with $k + \dim X \geq \ell$, and $X \cap F^k = \emptyset $ otherwise.

\begin{remark}
From now, we fix a generic flag $\cF$ sufficiently near to $\overline{H}_{\infty}$, that is,  $F^{k-1}$ does not separate $L_0 (\cA \cap F^{k})$ for each $0 \leq k \leq \ell$. 
\end{remark}

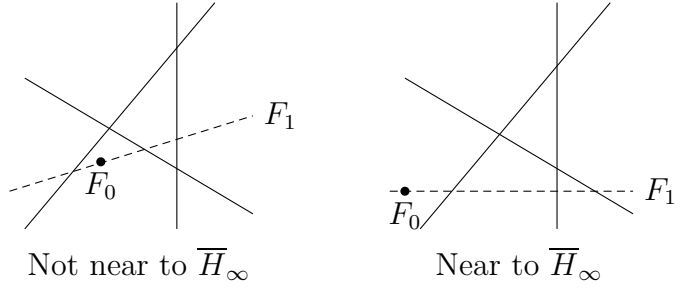
\begin{figure}[htbp]
\centering
\begin{tikzpicture}
\coordinate (A) at (0,0);
\coordinate (B) at (5,0);
\draw (A) ++ (0,0) --++ (2.5,3);
\draw (A) ++ (0,2) --++ (3,-1.8);
\draw (A) ++ (2,0) --++ (0,3);

\draw (B) ++ (0.2,0) --++ (2.5,3);
\draw (B) ++ (0,2) --++ (3,-1.8);
\draw (B) ++ (2,0) --++ (0,3);

\draw[densely dashed] (A) ++ (-0.2,0.5) --++ (3.2,1);
\draw[densely dashed] (B) ++ (-0.2,0.5) --++ (3.2,0);

\fill [black] (A)++(1,0.89) circle (0.06) node[below]{$F_0$};
\draw (A)++(3,1.5) node [right] {$F_1$} ;
\draw (A)++ (1.5,-0.1) node[below]{Not near to $\overline{H}_{\infty}$};

\fill [black] (B)++(0,0.5) circle (0.06) node[below]{$F_0$};
\draw (B)++(3,0.5) node [right] {$F_1$} ;
\draw (B)++ (1.5,-0.1) node[below]{Near to $\overline{H}_{\infty}$};

\end{tikzpicture}
\caption{Examples of flags}
\end{figure}

Define the subsets $\ch^{k} (\cA)$, $\bch^{k} (\cA)$, and $\uch^{k} (\cA)$ of $\ch(\cA)$ as follows:
\begin{eqnarray*}
\ch^{k} (\cA)& =& \{C \in \ch (\cA) \mid C \cap F^{k} \neq \emptyset, C \cap F^{k-1} = \emptyset\}, \\
\bch^{k} (\cA) &=& \{C\in \ch^{k} (\cA) \mid \mbox{$C \cap F^{k}$ is bounded} \}, \\
\uch^{k} (\cA) &=& \{C\in \ch^{k} (\cA) \mid \mbox{$C \cap F^{k}$ is unbounded} \}.
\end{eqnarray*}
From the assumption that $\cF$ is sufficiently near to $\overline{H}_{\infty}$, it follows that $\bch (\mathcal{A}) = \bch^{\ell} (\cA)$. 

It is known that the number of $\ch^{k} (\cA)$ is equal to the Betti number of the complexified complement.

\begin{proposition}\label{prop:betti} (\rm{\cite{yos-lef}, Proposition 2.3.2 (i)})
For $0 \leq k \leq \ell$, $b_k (M(\cA)) = \# \ch^{k} (\cA)$.
\end{proposition}

Let $C \in \uch({\cA})$. There exists a unique chamber $C^{\vee}$ which is the opposite with respect to $\overline{H}_{\infty}$. Since $(C^{\vee})^{\vee} = C$, we obtain a bijection 
\[
\iota : \uch (\cA) \rightarrow \uch (\cA), C \mapsto C^{\vee}.
\]

\begin{proposition} \label{prop:invol} \rm({\cite{yos-loc}, Theorem 2.10)} 
The involution $\iota : \uch (\cA) \rightarrow \uch (\cA)$ induces the bijection
\[
\iota : \bch^{k} (\cA) \xrightarrow{\cong} \uch^{k+1} (\cA) ,
\]
for $0 \leq k \leq \ell -1$. Thus it follows that $\# \bch^{k} (\cA) = \# \uch^{k+1} (\cA)$.
\end{proposition}

\begin{proposition}
Let $\beta_{k} = \# \bch^{k} (\cA)$. Then, $\beta_{k} = |\sum_{i=0}^{k} (-1)^{i} \cdot b_i (M(\cA))| $ for $0 \leq k \leq \ell-1$. In particular, $\# \bch(\cA)=\# \bch^{\ell} (\cA) =| \chi (M(\cA))|$.
\end{proposition}

\proof
It is clear that $\# \bch^{0} (\cA) = b_0 (M(\cA)) = 1$ for $k=0$.
When $k \geq 1$, from Proposition \ref{prop:betti} and \ref{prop:invol}, we obtain
\[
b_k (M(\cA)) = \# \ch^{k} (\cA) = \# \bch^{k} (\cA) + \# \uch^{k} (\cA) = \# \bch^{k} (\cA) + \# \bch^{k-1} (\cA).
\]
The proposition follows by induction on $k$.
\endproof

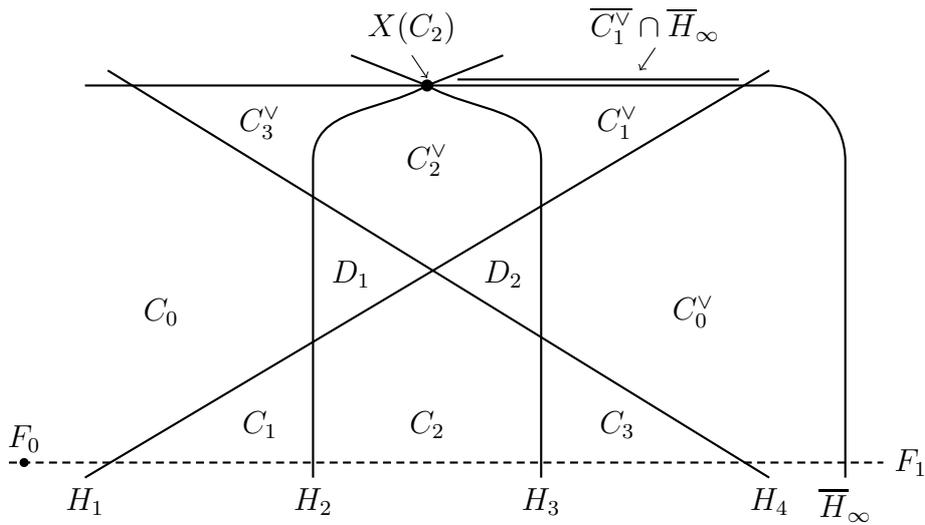
\begin{figure}[htbp]
\centering
\begin{tikzpicture}
\draw [thick, densely dashed](0,1)--(11.5,1);
\fill[black] (0.2,1) circle (0.06)  node[above] {$F_0$};
\draw (11.5,1) node[right]{$F_1$};

\draw [thick] (1,6) -- (10,6) to[out=0,in = 90] (11,5) -- (11,0.8) ;
\draw (11,0.8) node[below] {$\overline{H}_{\infty}$};

\draw [thick] (1,0.8)--(10,6.2);
\draw (1,0.8) node[below] {$H_1$};
\draw [thick] (4,0.8) -- (4,5)to [out=90, in=210] (5.5,6) to (6.5,6.4);
\draw (4,0.8) node[below] {$H_2$};
\draw [thick] (7,0.8) -- (7,5) to [out=90, in=330] (5.5,6) to (4.5,6.4);
\draw (7,0.8) node[below] {$H_3$};
\draw [thick] (1.3,6.2) -- (10,0.8);
\draw (10,0.8) node[below] {$H_4$};

\draw (2,3) node{$C_0$};
\draw (3.3,1.5) node{$C_1$};
\draw (5.5,1.5) node{$C_2$};
\draw (8,1.5) node{$C_3$};
\draw (9,3) node{$C_0^{\vee}$};
\draw (3.3,5.5) node{$C_3^{\vee}$};
\draw (5.5,5) node{$C_2^{\vee}$};
\draw (8,5.5) node{$C_1^{\vee}$};
\draw (4.5,3.5) node{$D_1$};
\draw (6.5,3.5) node{$D_2$};

\fill[black] (5.5,6) circle(0.08) ;
\draw [->] (5.3,6.4) -- (5.5,6.1);
\draw (5.3,6.4) node[above]{$X(C_2)$};

\draw [thick] (5.9,6.08) -- (9.6,6.08);
\draw [->] (8.5,6.5) -- (8.3,6.2);
\draw (8.5,6.4) node[above]{$\overline{C_1^{\vee}} \cap \overline{H}_{\infty}$};

\end{tikzpicture}
\caption{A flag and chambers}
\end{figure}

Next, we see some properties of unbounded chambers, their opposite chambers, and separating hyperplanes.
 We denote by $X(C)$ the minimum subspace containing $\overline{C} \cap \overline{H}_{\infty}$ for $C \in \uch (\cA)$.

\begin{proposition} \rm{(\cite{bai-yos}, Proposition 2.3)}
Let $C \in \uch (\cA)$ be an unbounded chamber. Then,
\[
\Sep (C,C^{\vee}) = \{H \in \cA \mid H \nsupseteq X(C)\} =\overline{\cA} \setminus \overline{\cA}_{X(C)}.
\]
In particular, if $\dim X(C) = \ell -1 $, then $\Sep (C,C^{\vee}) = \cA$.
\end{proposition}

\begin{proposition} \rm{(\cite{yos-loc}, Proposition 2.4)}
Let $X \in L(\overline{\cA})$ be an intersection contained in $\overline{H}_{\infty}$. Then, $X$ is a dense edge if and only if there exists a chamber $C \in \uch (\cA)$ such that $X = X(C)$. In particular, we have the following:
\[
D_{\infty} (\cA) = \{X(C) \mid C \in \uch (\cA)\}.
\]
\end{proposition}

From these propositions, we have the following equation of the monodromy.
\begin{proposition} \label{prop:monod}
Let $\cL$ be a rank one $R$-local system. For $C \in \uch (\cA)$, it follows that
\[
q_{\Sep (C,C^{\vee})} = q_{X(C)}^{-1}.
\]
In particular, the local system $\cL$ satisfies the CDO-condition if and only if $q_{\Sep (C,C^{\vee})} \neq 1$ for all $C \in \uch (\cA)$.
\end{proposition}

\subsection{Local system cohomology}

We now define the degree map $\deg : \ch^{k} (\cA) \times \ch^{k+1} (\cA) \rightarrow \bZ$. 
Let $C \in \ch^{k} (\cA)$ and $C' \in \ch^{k+1} (\cA)$. Let $B=B^{k} \subset F^{k}$ be a sufficiently large ball so that every $0$-dimensional edge $X \in L_0(\cA \cap F^{k})$ is contained in the interior of $B^{k}$. Then, there exists a tangent vector field $U^{C'}$ on $F^{k}$ that satisfies the following properties:
\begin{itemize}
\item[(i)] If $x \in \partial B \cap \overline{C}$, $U^{C'} (x) \neq 0$ and $U^{C'} (x)$ directs inside of $B$.
\item[(ii)] If $x \in H \cap F^{k}$ for $H \in \cA$, then $U^{C'} (x) \notin T_x (H \cap F^{k})$ and $U^{C'} (x)$ directs the side in which $C'$ is contained.  
\end{itemize}

\begin{definition}
Let $C\in \ch^{k} (\cA)$ and $C' \in \ch^{k+1}(\cA)$. Define the degree map $\deg (C, C')$ as follows.
\begin{itemize}
\item[(i)] When $k=0$, $\deg (C,C') = 1$.
\item[(ii)] When $k \geq 1$, 
\[
\deg (C,C') = \deg \left(\frac{U^{C'}}{|U^{C'}|} : \partial (\overline{C} \cap B) \rightarrow S^{k-1} \right) \in \bZ.
\]
\end{itemize}
\end{definition}

\begin{example}
Let us see the $k=1$ case. The degree map can be computed by using the combinatorial description of chambers.

\begin{itemize}
\item[(i)] Suppose $C \in \bch^{1} (\cA)$. Then $C \cap F^{1}$ is an open interval and the boundary is expressed as $(H \cap F^{1}) \cup (H' \cap F^{1})$ for $H , H' \in \cA$. Then $\deg (C, C')$ can be computed as
\begin{eqnarray*}
\deg(C,C') =
\left \{
\begin{array} {ll}
1& \mbox{if $H$, $H' \in \Sep (C,C')$}, \\
-1 & \mbox{if $H$, $H' \notin \Sep (C,C')$},  \\
0 & \mbox{otherwise}.
\end{array}
\right.
\end{eqnarray*}
\item[(ii)] Suppose $C \in \uch^{1}(\cA)$. Then $C \cap F^{1}$ is an unbounded interval and the boundary is a point which can be expressed as $H \cap F^{1}$ for $H \in \cA$. 
Then $\deg (C, C')$ can be computed as 
\begin{eqnarray*}
\deg(C,C') =
\left \{
\begin{array} {ll}
-1& \mbox{if $H \notin \Sep (C,C')$}, \\
0 & \mbox{if $H \in \Sep (C,C')$}.  
\end{array}
\right.
\end{eqnarray*}
\end{itemize}
\end{example}
\color{black}

It is difficult to calculate degrees in general. However, for the opposite chambers, the next proposition follows.
\begin{proposition} \rm {(\cite{bai-yos}, Theorem 4.8)} \label{prop:deg}
Let $C \in \bch^{k} (\cA)$ for $0 \leq k \leq \ell -1$. Then, 
\[
\deg (C,C^{\vee}) = (-1)^{\ell - 1 -\dim X(C)} .
\]
\end{proposition}

Let $R[\ch^{k} (\cA)] = \oplus_{C \in \ch^{k} (\cA)} R[C]$ be the free $R$-module with the basis $\ch^{k} (\cA)$, where $R$ is the coefficient ring.
For the chambers $C$, $C' \in \ch(\cA)$ define $\Delta (C,C') \in R$ as follows:
\begin{equation}\label{def:delta}
\Delta (C,C') = \sqrt{-1} \cdot (q^{1/2}_{\Sep (C,C')} - q^{-1/2}_{\Sep (C,C')}),
\end{equation}
where $q^{1/2}_{\Sep (C,C')} = \prod_{H_i \in \Sep (C,C')} q_i ^{1/2}$. 
Since $q_i = \pm1 $ for $\bZ$-local systems, it is clear that $q_i^{1/2} - q_i^{-1/2} =\pm 2  \cdot \sqrt{-1}$ or $0$.
Therefore, $\Delta (C,C') = \pm 2$ or $0$. Also note that $\Delta (C,C') = 0$ if and only if $q_{\Sep (C,C')} = 1$.

\begin{remark}
Let $t(C,C')$ be the number of “twisted" hyperplanes separating $C$ and $C'$, that is,
\[
t (C,C') = \# \{q_i \mid \mbox{$H_i \in \Sep (C,C')$ , $q_i = - 1$}\}.
\]
For all $i$ with $q_i=-1$, let us fix the branch $q_i^{1/2} = \sqrt{-1}$.
Then, we can determine $\Delta (C, C')$ for $\bZ$-local system by counting “twisted" hyperplanes mod 4:
\begin{eqnarray*}
\Delta (C,C') = \left \{ 
\begin{array} {ll}
-2 & (\mbox{$t(C,C') \equiv 1$ mod $4$}),\\
2 & (\mbox{$t(C,C') \equiv 3$ mod $4$}),\\
0 & (\mbox{$t(C,C') \equiv 0,2$ mod $4$}).
\end{array}
\right.
\end{eqnarray*}

\end{remark}

\begin{definition}
Define the $R$-homomorphism $d_{\cL} ; R[\ch^{q}(\cA)] \rightarrow R[\ch^{q+1} (\cA)]$ by
\[
d_{\cL} ([C]) = \sum_{C' \in \ch^{q+1} (\cA)} \deg (C,C') \cdot \Delta (C,C') \cdot [C'].
\]
\end{definition}

\begin{proposition} \rm{(\cite{yos-lef}, Theorem 6.4.1 for $R=\bC$ case)} 
$(R[\ch^{\bullet} (\cA)],d_{\cL})$ is a cochain complex. Furthermore,
\[
H^{*} (R[\ch^{\bullet} (\cA)], d_{\cL}) \cong  H^{*} (M(\cA), \cL) .
\]
\end{proposition}

\begin{remark}
In \cite{yos-lef}, only the computation for $\bC$-local system cohomology is described. Actually, we can see that Theorem 6.4.1 in \cite{yos-lef} holds also for $\bZ$-local system from a similar argument. However, only the definition of $\Delta_{\Phi} (C, C')$ in that paper looks different from the one defined in (\ref{def:delta}) this paper. 
From the following discussion on the representation of a Deligne groupoid, it follows that these definitions are consistent.
Most notations are from Sections 3.3 and 6.5 in \cite{yos-lef}.

In \cite{yos-lef}, an explicit cell decomposition of $M(\cA)$ is obtained by using the combinatorial structure of chambers. 
The attaching map of a cell corresponding to a chamber $C$ is denoted by $\sigma_{C}$ (see Section 5 in \cite{yos-lef} for its definition and construction).
Let $\Gal (\cA)$ be a Deligne groupoid and $C,C' \in \ch (\cA)$.
The equivalence class of gallery defined by a geodesic from $C$ to $C'$ is denoted by $P^{+} (C,C') \in \Hom_{\Gal} (C,C')$. The inverse is denoted by $P^{-} (C',C) = (P^{+} (C,C'))^{-1} \in \Hom_{\Gal} (C',C)$. 

Let us fix numbers $q^{+}_{1} , \cdots, q^{+}_{n}, q^{-}_{1}, \cdots , q^{-}_{n} \in \{\pm 1\}$.
Then, we can define a representation $\Phi: \Gal (\cA) \rightarrow \bZ \mbox{-} mod$ as follows.
First, for each chamber, we define
\[
\Phi (C) = \bZ [\sigma_{C}],
\]
where $\bZ[\sigma_{C}]$ is a rank one free $\bZ$-module generated by $[\sigma_{C}]$. 

We suppose that $H_i$ is defined by a linear equation $\{ \alpha_{i}=0 \}$ and $\Sep (C,C') = \{H_{i_{1}} , \cdots , H_{i_{k}}\}$ for $C,C' \in \ch (\cA)$.
Then, define $\Phi_{P^{+} (C,C')}$ as
\[
\Phi_{P^{+}(C,C')} : \Phi (C) \rightarrow \Phi (C'), [\sigma_{C}] \mapsto q^{\varepsilon_{i_{1}}}_{i_{1}} q^{\varepsilon_{i_{2}}}_{i_{2}} \cdots q^{\varepsilon_{i_{k}}}_{i_{k}} [\sigma_{C'}],
\]
where $\varepsilon_{i_{j}} = +$ ($-$) if a chamber moves from the side $\{\alpha_{i_{j}} <0 \mbox{ ($>0$)}\}$ to the side $\{\alpha_{i_{j}} >0 \mbox{ ($<0$)}\}$.
A morphism $P^{+}(C',C)$ is represented by a gallery of the reverse orientation of that of $P^{+} (C,C')$. 
Therefore, $\Phi_{P^{+} (C',C)}$ is written as 
\[
\Phi_{P^{+}(C',C)} : \Phi (C') \rightarrow \Phi (C), [\sigma_{C'}] \mapsto q^{-\varepsilon_{i_{k}}}_{i_{k}} q^{-\varepsilon_{i_{k-1}}}_{i_{k-1}} \cdots q^{-\varepsilon_{i_{1}}}_{i_{1}} [\sigma_{C}].
\]
Since $P^{-} (C,C') = (P^{+} (C',C))^{-1}$ and the functoriality, $\Phi_{P^{+} (C',C)}$ can be written as
\begin{eqnarray*}
\Phi_{P^{-}(C,C')} : \Phi (C) \rightarrow \Phi (C'), [\sigma_{C}] &\mapsto& \left( q^{-\varepsilon_{i_{1}}}_{i_{1}} \right)^{-1} \left(q^{-\varepsilon_{i_{2}}}_{i_{2}} \right)^{-1} \cdots \left(q^{-\varepsilon_{i_{k}}}_{i_{k}} \right)^{-1} [\sigma_{C'}] \\
&=& q^{-\varepsilon_{i_{1}}}_{i_{1}}  q^{-\varepsilon_{i_{2}}}_{i_{2}} \cdots q^{-\varepsilon_{i_{k}}}_{i_{k}} [\sigma_{C'}].
\end{eqnarray*}


The represetation defined in this way determines a rank one $\bZ$-local system such that local monodromy around $H_i$ is $q^{+}_{i} q^{-}_{i}$. 
Note that local monodromy $q_i = 1$ if and only if $q^{+}_{i} = q^{-}_{i}$. 
Under this setting, $\Delta_{\Phi} (C,C') = \Phi_{P^{+}(C,C')} - \Phi_{P^{-}(C,C')}$ is written as
\[
\Delta_{\Phi} (C,C') = q^{\varepsilon_{i_{1}}}_{i_{1}} q^{\varepsilon_{i_{2}}}_{i_{2}} \cdots q^{\varepsilon_{i_{k}}}_{i_{k}} -q^{-\varepsilon_{i_{1}}}_{i_1} q^{-\varepsilon_{i_{2}}}_{i_{2}} \cdots q^{-\varepsilon_{i_{k}} }_{i_{k}}.
\]
It is easy to see that this coincides with the one defined in (\ref{def:delta}) this paper. 
In fact, $\Delta_{\Phi} (C,C') = \pm 2$ or $0$ and $\Delta_{\Phi} (C,C') = 0$ if and only if $q_{\Sep (C,C')} = 1$.
\end{remark}

We now give an example of the computation of $\bZ$-local system cohomology.

\begin{example}
Let $\cA = \{H_1, H_2, H_3\}$ be the arrangement of central $3$ lines in $\bR^2$. Suppose that generic flag $F^{0} \subset F^{1}$ is located as in the figure.

\begin{figure}[htbp]
\centering
\begin{tikzpicture}
\draw [thick](-3,0)--(3,4);
\draw (-3,0) node[below] {$H_1$};

\draw [thick](0,0)--(0,4);
\draw (0,0) node[below] {$H_2$};

\draw [thick](-3,4) -- (3,0);
\draw (3,0) node[below] {$H_3$};

\draw [thick, densely dashed] (-4,0.3) --(3.5,0.3);
\draw (3.5,0.3) node[right] {$F_1$};
\fill[black] (-3.5,0.3) circle (0.06) node[above] {$F_0$};

\draw (-2,2) node{$C_0$};
\draw (-0.8,1) node{$C_1$};
\draw (0.8,1) node{$C_2$};
\draw (2,2) node{$C_0^{\vee}$};
\draw (-0.8,3.5) node{$C_2^{\vee}$};
\draw (0.8,3.5) node{$C_1^{\vee}$};

\end{tikzpicture}
\caption{Central $3$ lines}
\end{figure}
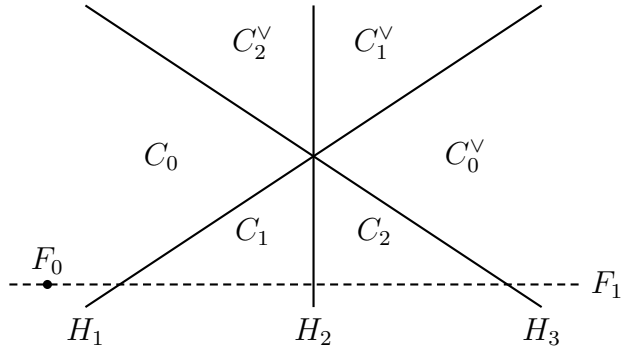

The boundary map of the cochain complex $(\bZ[\ch^{\bullet}(\cA)], d_{\cL})$ is expressed as follows:
\begin{eqnarray*}
d_{\cL} ([C_0]) &=& \sqrt{-1} \cdot ( (q_1^{1/2} - q_1^{-1/2} )[C_1] + (q_{12}^{1/2} - q_{12}^{-1/2} ) [C_2] + (q_{123}^{1/2} - q_{123}^{-1/2} )[C_0^{\vee}] ) \\ 
d_{\cL} ([C_1] ) &=& \sqrt{-1} \cdot (q_{123}^{1/2} - q_{123}^{-1/2}) [C_1 ^ {\vee}] \\
d_{\cL} ([C_2] ) &=& \sqrt{-1} \cdot (q_{123}^{1/2} - q_{123}^{-1/2}) [C_2 ^ {\vee}] \\
d_{\cL} ([C_0 ^{\vee}] ) &=&  \sqrt{-1} \cdot (-(q_{1}^{1/2} - q_{1}^{-1/2}) [C_1 ^ {\vee}] -(q_{12}^{1/2} - q_{12}^{-1/2}) [C_2 ^ {\vee}] )
\end{eqnarray*}
By computing them, we have the following:

\begin{table}[htbp]
\centering
\begin{tabular}{|c||c|c|c|} \hline
&$q_i = 1, \forall i$&$q_1 q_2 q_3 = 1, \mbox{$\exists i$ s.t. $q_i \neq 1$} $&otherwise \\ \hline
$H^0 (\bZ[\ch^{\bullet} (\cA)],d_{\cL})$ & $\bZ$ & $0$ & $0$\\  \hline
$H^1 (\bZ[\ch^{\bullet} (\cA)],d_{\cL})$ & $\bZ^3$ & $\bZ \oplus \bZ / 2\bZ$ & $\bZ / 2\bZ$\\  \hline
$H^2 (\bZ[\ch^{\bullet} (\cA)],d_{\cL})$ & $\bZ^2$ & $\bZ \oplus \bZ / 2\bZ$ & $\bZ / 2\bZ \oplus \bZ / 2\bZ$\\  \hline

\end{tabular}
\caption{Local system cohomology via chamber basis}
\end{table}

On the other hand, the complexified complement $M(\cA)$ has the homotopy type of $2$-dimensional CW complex obtained by the following presentation of the fundamental group. (\cite{coh-suc-braid} \cite{yos-lef})
\[
\pi_1 (M(\cA)) \cong \langle \gamma_1, \gamma_2, \gamma_3 \mid \gamma_1 \gamma_2 \gamma_3 \gamma_1^{-1} \gamma_3^{-1} \gamma_2^{-1}, \gamma_1 \gamma_2 \gamma_3 \gamma_2^{-1} \gamma_1^{-1} \gamma_3^{-1}   \rangle
\]
We can compute local system cohomology directly from this presentation and have consistent results.
\end{example}

\section{The main result and its proof}




The local system $\cL$ is always a rank one $\bZ$-local system throughout this section.
We give the proof of our main result in this section. 
First, let us repeat the main result. 

\begin{theorem} \label{thm:main}
Let $\cA = \{H_1 , \cdots , H_n\}$ be a real hyperplane arrangement in $\bR^{\ell}$ and denote by $M(\cA) $ the complexified complement. Suppose that the $\bZ$-local system $\cL$ satisfies CDO-condition. Then,
\begin{eqnarray*}
H^{k} (M(\cA), \cL) \cong
\left \{
\begin{array} {ll}
\bZ^{\beta_{\ell}} \oplus \bZ_2^{\beta_{\ell -1}} & (k= \ell ), \\
\bZ_2^{\beta_{k-1}} & (1 \leq k \leq \ell-1), \\
0 & (otherwise),
\end{array}
\right.
\end{eqnarray*}
where $\beta_{k} = \# \bch^{k} (\cA) = |\sum_{i=0}^{k} (-1)^i \cdot b_i (M(\cA))|$.
\end{theorem}

Define the homomorphism $d'_{\cL} : \bZ[\ch^{k} (\cA)] \rightarrow \bZ[\ch^{k+1} (\cA)]$ as follows:
\[
d'_{\cL} (C) =\sum_{C' \in \ch^{k+1} (\cA)} \frac{1}{2} \cdot \deg(C,C') \cdot \Delta (C,C') \cdot [C'].
\]
Since $\Delta (C,C') =\pm 2$ or $0$, this is a well-defined homomorphism. It is easy to check that $d'^{k+1}_{\cL} \circ d'^{k}_{\cL} = 0$ and $\Ker (d'^{k}_{\cL}) = \Ker (d^{k}_{\cL})$. Next, we define the homomorphism $\overline{d_{\cL}} $ as the composition
\[
\overline{d_{\cL}} : \bZ[\bch^{k} (\cA)] \hookrightarrow \bZ[\ch^{k} (\cA)] \xrightarrow{d_\cL} \bZ[\ch^{k+1} (\cA) ] \twoheadrightarrow \bZ[\uch^{k+1} (\cA)].
\]
Since $\# \bch^{k} (\cA) = \# \uch^{k+1} (\cA)$ for each $k$ (Proposition \ref{prop:invol}), $\overline{d_{\cL}}$ is represented by a square matrix by fixing the basis $\bch^{k} (\cA)= \{C_1 , \cdots, C_b\}$. 
We define $\overline{d'_{\cL}}$ for $d'_{\cL}$ similarly. The next proposition implies that the homomorphisms $\overline{d_{\cL}}$ and $\overline{d'_{\cL}}$ can be represented by an upper triangular matrix by reordering the chamber basis appropriately.

\begin{proposition} \rm{(\cite{bai-yos} Theorem 4.6)} \label{prop:upper}
For $0 \leq k \leq \ell-1$, let us fix an ordering of chambers $\bch^{k}(\cA) = \{C_1 , \cdots , C_b\}$ so that
\[
\dim X(C_1) \geq \dim X(C_2) \geq \cdots \geq \dim X(C_b).
\]
Then, the matrix $(\deg (C_i,C_j^{\vee}) )_{i,j}$ is an upper triangular, that is, if $i > j$ , then $\deg (C_i, C_j^{\vee}) = 0$.
\end{proposition}

The CDO-condition is crucial for the non-degeneracy of the homomorphism $\overline{d_{\cL}}$ and $\overline{d'_{\cL}}$.

\begin{proposition}
Suppose that local system $\cL$ satisfies CDO-condition. Then, the homomorphisms $\overline{d_{\cL}}$, $\overline{d'_{\cL}}$ are both injective. Furthermore, $\overline{d'_{\cL}}$ is an isomorphism.
\end{proposition}

\proof
Let us fix the ordering of chambers in $\bch^{k}(\cA)$ as in Proposition \ref{prop:upper} . Then, since the coefficient matrix is upper triangular, the determinant of the homomorphism is the product of all diagonal elements. 
Each diagonal element of $\overline{d_{\cL}}$ is expressed as
\[
(\overline{d_{\cL}})_{(i,i)} = \sqrt{-1} \cdot \deg (C_i,C_i^{\vee})\cdot (\pm (q_{\Sep(C_i,C_i^{\vee})}^{1/2} - q_{\Sep(C_i,C_i^{\vee})}^{-1/2})).
\]
By CDO-condition, we have $q_{\Sep(C_i,C_i^{\vee})}^{1/2} - q_{\Sep(C_i,C_i^{\vee})}^{-1/2} = \pm 2 \sqrt{-1} \neq 0$ for each $C_{i} \in \bch^{k} (\cA)$. Therefore the product of diagonal elements is non zero (see Proposition \ref{prop:monod} and  \ref{prop:deg}). This implies the injectivity of $\overline{d_{\cL}}$. Since each diagonal element of $\overline{d'_{\cL}}$ is $\pm 1$, the determinant $\det (\overline{d'_{\cL}}) = \pm 1$. Therefore, $\overline{d'_{\cL}}$ defines an isomorphism.
\endproof

\begin{proposition} \label{prop:basis}
For $1 \leq k \leq \ell$, the set $\bch^{k} (\cA) \cup d'_{\cL} (\bch^{k-1} (\cA))$ forms a basis of $\bZ [\ch^{k} (\cA)]$.
\end{proposition}



\proof 
We give the proof by induction on $k$. When $k=1$, let $\ch^{0} (\cA) = \{C_0\}$ and $\bch^{1} (\cA) = \{C_1, \cdots ,C_{n-1}\}$. The image $d'_{\cL} ([C_0])$ is expressed as
\[
d'_{\cL} ([C_0]) = \frac{\sqrt{-1}}{2} \cdot \left( (q_1^{1/2} - q_1^{-1/2}) [C_1] + \cdots + (- (q_{\infty}^{1/2} - q_{\infty}^{-1/2}) )[C_0^{\vee}] \right).
\]
By CDO condition, the coefficient of $[C_0^{\vee}]$ is $\pm 1$. Therefore, it is clear that $\{C_1 , \cdots, C_{n-1} , d'_{\cL} ([C_0])\}$ forms a basis of $\bZ[\ch^{1} (\cA)]$. 

Suppose $k > 1$ and let $\bch^{k-1} (\cA) = \{C_1 , \cdots, C_s\}$ and $\bch^{k} (\cA) = \{D_1 , \cdots, D_t\}$. 
By the definition of the homomorphism $\overline{d'_{\cL}}$, there exists an element $A_i \in \bZ[\bch^{k}(\cA)]$ such that $d'_{\cL} ([C_i] ) = A_i + \overline {d'_{\cL}} ([C_i])$. 
First, we prove linear independence.
Consider the linear equation
\[
\sum_{i=1}^{s} a_i \cdot d'_\cL ([C_i] )+ \sum_{j=1}^{t} b_j \cdot ( [D_j] ) = 0,
\]
where $a_i$, $b_j$ are integers. This is equivalent to 
\[
\sum_{i=1}^{s} a_i \cdot \left( A_i + \overline{d'_\cL} ([C_i]) \right) + \sum_{j=1}^{t} b_j \cdot ([D_j]) = 0. 
\]
The elements in $\bZ[\uch^{k}(\cA)]$ appear only in the terms of $\overline{d'_{\cL}} ([C_i])$'s. Therefore, it follows that $a_i = 0$ for all $i$ from the injectivity of $\overline{d'_{\cL}}$. Then all $b_j$ must be $0$. 
This completes the proof of independence. 
To prove that the set $\bch^{k} (\cA) \cup d'_{\cL} (\bch^{k-1} (\cA) )$ forms a generating set of $\bZ[\ch^{k} (\cA)]$, it is enough to show that $C_j^{\vee} \in \bZ[\uch^{k} (\cA)]$ is generated by the set. 
Since $\overline{d'_{\cL}} : \bZ[\bch^{k-1} (\cA)] \rightarrow \bZ[\uch^{k}(\cA)]$ is isomorphism, there exist unique integers $a_1, \cdots , a_n \in \bZ$ such that 
\[
[C^{\vee}_j] = \sum_{i=1}^{s} a_i \cdot \overline{d'_{\cL}}([C_i]).
\]
Since $\overline{d'_{\cL}} ([C_i]) = d'_{\cL} ([C_i]) - A_i$ and $A_i \in \bZ[\bch^{k} (\cA)]$, $[C^{\vee}_i]$ can be written by the linear combination of the elements in the set $\bch^{k} (\cA) \cup d'_{\cL} (\bch^{k-1} (\cA))$.

\endproof

\begin{proposition} \label{prop:isom}
For $1 \leq k \leq \ell$, $\Ker (d'^{k}_{\cL}) = \Img (d'^{k-1}_{\cL}) \cong \bZ^{\beta_{k-1}}$.
\end{proposition}

\proof
Let $\bch^{k-1}(\cA) = \{C_1 ,\cdots , C_s\}$ and $\bch^{k} (\cA) = \{D_1 ,\cdots D_t\}$.
It is clear that $\Img (d'_{\cL}) \subset \Ker (d'_{\cL})$. Therefore we need to prove the inverse inclusion. 
From Proposition \ref{prop:basis}, an arbitrary element $A \in \bZ[\ch^{k} (\cA)]$ is written as $A = \sum_{i=1}^{s} a_i \cdot d'_{\cL} ([C_i]) + \sum_{j=1}^{t} b_j \cdot [D_j ]$ uniquely ($a_i, b_j \in \bZ$). 
Since $d'_{\cL} \circ d'_{\cL} = 0$, the image $d'_{\cL} (A)$ is written as follows :
\[
d'_{\cL} (A) = \sum_{j=1}^{t} b_j \cdot d'_{\cL} ([D_j]).
\]
If we assume that $A \in \Ker (d'_{\cL})$, then $b_j = 0$ for all $j$ from the injectivity of $\overline{d'_{\cL}}$. Therefore, any element of $\Ker (d'_{\cL})$ is written as $A = \sum_{i=1}^{s} a_i \cdot d'_{\cL} ([C_i])$. This is contained in $\Img (d'_{\cL})$. Finally, the proposition follows from that $\Img (d'_{\cL}) = \oplus_{C \in \bch^{k-1} (\cA)} \bZ (d'_{\cL}([C])) \cong \bZ^{\beta_{k-1}}$.

\endproof

\noindent
{\bf The proof of Theorem \ref{thm:main}.} 
\begin{itemize}
\item[(i)]When $k=0$, $d_\cL : \bZ[\ch^{0} (\cA)] \rightarrow \bZ[\ch^{1} (\cA)]$ is injective. Therefore $H^{0} (\bZ[\ch^{\bullet} (\cA)] , d_{\cL}) = 0$.
\item[(ii)]  When $1 \leq k \leq \ell-1$, $\Ker (d'_{\cL}) = \Ker (d_{\cL})$, $d_{\cL} = 2 \cdot d'_{\cL}$, and Proposition \ref{prop:isom} implies
\[
H^k(\bZ[\ch^{\bullet}(\cA)],d_{\cL}) \cong \frac{\Ker(d^k_{\cL})}{\Img (d^{k-1}_{\cL})} \cong \frac {\Ker(d'^k_{\cL})}{2 \Img (d'^{k-1}_{\cL})} \cong \frac {\Ker(d'^k_{\cL})}{2\Ker(d'^k_{\cL} )} \cong \bZ_2^{\beta_{k-1}}.
\]
\item[(iii)] When $k=\ell$, we obtain the following similarly as (ii):
\[
H^{\ell} (\bZ[\ch^{\bullet} (\cA),d_{\cL}]) \cong \frac{\bZ[\ch^{\ell} (\cA)]}{\Img(d^{\ell-1}_{\cL})} \cong \frac{\bZ[\bch^{\ell} (\cA)] \oplus \Img (d'_{\cL})}{2 \Img (d'_{\cL})} \cong \bZ^{\beta_{\ell}} \oplus \bZ_2 ^{\beta_{\ell-1}} .
\]
\end{itemize}

\qed



\section{Examples --$2$-torsions in double coverings--}
Let $\omega : \pi_1 (M(\cA)) \rightarrow \bZ_2$ be a (surjective) group homomorphism. Then, $\omega$ defines a double covering $p(\omega) : M(\cA)^{\omega} \rightarrow M(\cA)$. 
Since there is a natural isomorphism
\[
\Hom (\pi_1 (M(\cA)), \bZ_2) \cong \Hom (H_1 (M(\cA)), \bZ_2) \cong H^1 (M(\cA), \bZ_2),
\]
this is equivalent to determine a (nonzero) element $\omega \in H^{1} (M(\cA),\bZ_2)$. This is called the characteristic class of the double covering. 
Since $\bZ_2 \cong \bZ^{\times}$, we can define the associated rank one $\bZ$-local system $\cL_{\omega}$ from the double covering.

The icosidodecahedral arrangement is a central hyperplane arrangement $\cA_{\cI \cD} = \{H_1, \cdots , H_{16}\}$ in $\bC^3$ arising from an icosidodecahedron which is introduced by Yoshinaga (\cite{yos-dou}). The icosidodecahedral arrangement is known to have the following two properties:
\begin{itemize}
\item[$\bullet$] The first integral homology of the Milnor fiber $H_{1} (F_{\cA_{\cI \cD}},\bZ)$ has a non-trivial $2$-torsion. 
\item[$\bullet$] A counterexample to the formula for the dimension of $(-1)$-eigenspace of Milnor fiber homology conjectured by Papadima-Suciu. (Conjecture 1.9 in \cite{pap-suc}).
\end{itemize} 



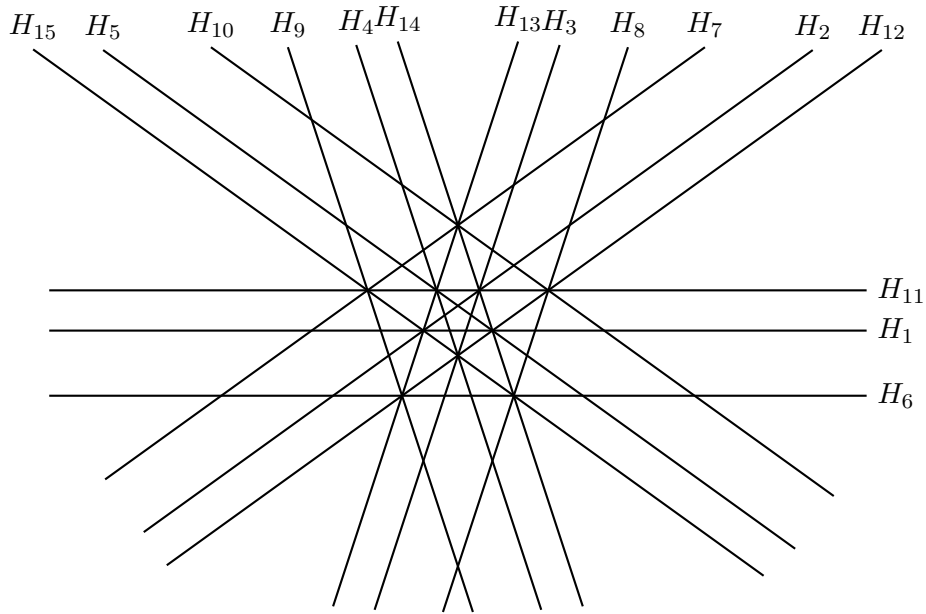
\begin{figure}[htbp]
\centering
\begin{tikzpicture}[scale=1.25]
\small




\draw[thick,rotate=-18] (0.809,0)++(0,-3) -- +(90:6.3) node[above]{$H_8$};
\draw[thick,rotate=-18] (-0.309,0)++(0,-3.3) -- +(90:6.3) node[above]{$H_{13}$};
\draw[thick,rotate=-18] (0.118,0)++(0,-3.2) -- +(90:6.3) node[above]{$H_3$};

\draw[thick,rotate=-18] (72:0.809)++(162:3.8)  node[above]{$H_{10}$}-- +(342:8.1);
\draw[thick,rotate=-18] (252:0.309)++(162:5.3)  node[above]{$H_{15}$}-- +(342:9.5);
\draw[thick,rotate=-18] (72:0.118)++(162:4.7)  node[above]{$H_5$}-- +(342:9);

\draw[thick,rotate=-18] (144:0.809)++(234:4) -- +(54:7.8) node[above]{$H_7$};
\draw[thick,rotate=-18] (322:0.309)++(234:4) -- +(54:9.3) node[above]{$H_{12}$};
\draw[thick,rotate=-18] (144:0.118)++(234:4) -- +(54:8.7) node[above]{$H_2$};

\draw[thick,rotate=-18] (216:0.809)++(126:3.3)  node[above]{$H_9$}-- +(306:6.3);
\draw[thick,rotate=-18] (36:0.309)++(126:3)  node[above]{$H_{14}$}-- +(306:6.3);
\draw[thick,rotate=-18] (216:0.118)++(126:3.1)  node[above]{$H_4$}-- +(306:6.3);

\draw[thick,rotate=-18] (288:0.809)++(198:4.3) -- +(18:8.6) node[right]{$H_6$};
\draw[thick,rotate=-18] (108:0.309)++(198:4.3) -- +(18:8.6) node[right]{$H_{11}$};
\draw[thick,rotate=-18] (288:0.118)++(198:4.3) -- +(18:8.6) node[right]{$H_1$};

\end{tikzpicture}
\caption{A deconing of the icosidodecahedral arrangement 
$\cA_{\cI \cD}=\{H_1, \dots, H_{15}\}$} 
\label{fig:ID}
\end{figure}

Torsions in the first homology groups Milnor fiber of multiarrangements are studied in \cite{cds,den-suc-multi,den-suc-tor}, however, the icosidodecahedral arrangement is the first example whose first homology of the Milnor fiber has a non-trivial $2$-torsion.

Recently, it is revealed that the two properties above are closely related (\cite{ish-sug-yos}). 
(A deconing of) the icosidodecahedral arrangement $\cA_{\cI \cD}$ also gives the example such that the first homology of double covering has a $2$-torsion. Moreover, $\bZ_{4}$-summand appears in the associated rank one $\bZ$-local system homology, which is never seen if the local system satisfies the CDO-condition.

\begin{example} {\rm (\cite{yos-dou,ish-sug-yos})}
Let $\cA_{\cI \cD} = \{H_1 , \cdots , H_{15}\}$ be a deconing of the icosidodecaheral arrangement. Then, there are some choices of characteristic classes $\omega \in H^1 (M(\cA_{\cI \cD}), \bZ_2 )$ such that the first homology of associated double covering has a non-trivial torsion and the associated local system homology has a $\bZ_4$-summand simultaneously.  Precisely, the following holds for such characteristic classes:
\begin{eqnarray*}
H_{1} (M(\cA_{\cI \cD})^{\omega}, \bZ) \cong \bZ^{15} \oplus \bZ_{2}, \\
H_1 (M(\cA_{\cI \cD}), \cL_{\omega}) \cong \bZ_{2}^{13} \oplus \bZ_{4}.
\end{eqnarray*}

\end{example}

In \cite{ish-sug-yos}, a smaller example with similar properties was found.




\begin{figure}[htbp]
\centering
\begin{tikzpicture}[scale=1.25]
\small




\draw[thick,rotate=-18] (0.809,0)++(0,-3) -- +(90:6.3) node[above]{$H_8$};
\draw[thick,rotate=-18] (-0.309,0)++(0,-3.3) -- +(90:6.3) node[above]{$H_{13}$};

\draw[thick,rotate=-18] (72:0.809)++(162:3.8)  node[above]{$H_{10}$}-- +(342:8.1);
\draw[thick,rotate=-18] (252:0.309)++(162:5.3)  node[above]{$H_{15}$}-- +(342:9.5);

\draw[thick,rotate=-18] (144:0.809)++(234:4) -- +(54:7.8) node[above]{$H_7$};
\draw[thick,rotate=-18] (322:0.309)++(234:4) -- +(54:9.3) node[above]{$H_{12}$};

\draw[thick,rotate=-18] (216:0.809)++(126:3.3)  node[above]{$H_9$}-- +(306:6.3);
\draw[thick,rotate=-18] (36:0.309)++(126:3)  node[above]{$H_{14}$}-- +(306:6.3);

\draw[thick,rotate=-18] (288:0.809)++(198:4.3) -- +(18:8.6) node[right]{$H_6$};
\draw[thick,rotate=-18] (108:0.309)++(198:4.3) -- +(18:8.6) node[right]{$H_{11}$};

\end{tikzpicture}
\caption{The double star arrangement
$\cA_{\cD \cS}=\{H_6, \cdots, H_{15}\}$} 
\label{fig:DP}
\end{figure}

\begin{example} {\rm (\cite{ish-sug-yos})}
Let $\cA_{\cD \cS} = \{H_6, \cdots, H_{15}\}$ be the subarragement of $\cA_{\cI \cD}$ as in the Figure \ref{fig:DP}. There is just one characteristic class $\omega \in H^{1} (M(\cA_{\cD \cS}),\bZ_2)$ satisfying the following:
\begin{eqnarray*}
H_{1} (M(\cA_{\cD \cS})^{\omega}, \bZ) \cong \bZ^{10} \oplus \bZ_{2}, \\
H_1 (M(\cA_{\cD \cS}), \cL_{\omega}) \cong \bZ_{2}^{8} \oplus \bZ_{4}.
\end{eqnarray*}
\end{example}

In the examples known to the author, non-trivial $2$-torsion of the first homology of the double covering appears if and only if the associated rank one $\bZ$-local system homology has a $\bZ_4$-summand.


\end{document}